\documentclass[12pt]{article}
\usepackage{amsthm,amsfonts,amssymb,amscd,amsmath}

\begin{document}

\begin{center}
\large \bf Effective birational rigidity of \\ Fano double
hypersurfaces
\end{center}\vspace{0.5cm}

\centerline{Thomas Eckl and Aleksandr Pukhlikov}\vspace{0.5cm}

\parshape=1
3cm 10cm \noindent {\small \quad\quad\quad \quad\quad\quad\quad
\quad\quad\quad {\bf }\newline We prove birational superrigidity
of Fano double hypersurfaces of index one with quadratic and
multi-quadratic singularities, satisfying certain regularity
conditions, and give an effective explicit lower bound for the
codimension of the set of non-rigid varieties in the natural
parameter space of the family. The lower bound is quadratic in the
dimension of the variety. The proof is based on the techniques of
hypertangent divisors combined with the recently discovered
$4n^2$-inequality for complete intersection
singularities.\vspace{0.1cm}

Bibliography: 18 titles.} \vspace{1cm}

\vspace{1cm}

\noindent 14E05, 14E07\vspace{1cm}

\section*{Introduction}

{\bf 0.1. Statement of the main result.} Fix the integers
$M\geqslant 10$, $m\geqslant 2$ and $l\geqslant 2$, satisfying the
equality
$$
m+l=M+1.
$$
Let ${\mathbb P}={\mathbb P}^{M+1}$ be the complex projective
space. By the symbol ${\cal P}_{k,M+2}$ we denote the space of
homogeneous polynomials of degree $k\in {\mathbb Z}_+$ in $M+2$
homogeneous coordinates on ${\mathbb P}$, that is, the linear
space $H^0({\mathbb P},{\cal O}_{\mathbb P}(k))$. Let
$$
(g,h)\in {\cal P}_{m,M+2}\times {\cal P}_{2l,M+2}={\cal F}
$$
be a pair of irreducible polynomials.\vspace{0.1cm}

Consider the double cover
$$
\sigma\colon V\stackrel{2:1}{\longrightarrow} G\subset {\mathbb
P},
$$
where $G=\{g=0\}\subset {\mathbb P}$ is an irreducible
hypersurface of degree $m$ and $\sigma$ is branched over the
divisor $W=\left\{h|_G=0\right\}\subset G$, which is cut out on
$G$ by the hypersurface $W_{\mathbb P}=\{h=0\}$. The variety $V$
can be realized as a complete intersection of codimension 2 in the
weighted projective space
$$
{\mathbb P}(\underbrace{1,\dots,1}_{M+2}, l)={\mathbb
P}(1^{M+2},l)
$$
with the homogeneous coordinates $x_0,\dots,x_{M+1}$ of weight 1
and the new homogeneous coordinate $u$ of weight $l$:
$$
V=\{g=0,\,\, u^2=h\}.
$$
If the variety $V$ is factorial and its singularities are
terminal, then $V$ is a primitive Fano variety:
$$
\mathop{\rm Pic} V={\mathbb Z}H,\quad K_V=-H,
$$
where $H$ is the class of ``hyperplane section'', corresponding to
$\sigma^* {\cal O}_{\mathbb P}(1)|_G$. It makes sense now to test
$V$ for being birationally (super)rigid. In \cite{Pukh00e} it was
shown that a Zariski general non-singular variety $V$ is
birationally superrigid. The aim of this paper is to generalize
and strengthen that result in the following way.\vspace{0.1cm}

Let us define the integer-valued function
$$
\xi\colon {\mathbb Z}_{\geqslant 10} =\{M\in {\mathbb Z}\,|\,
M\geqslant 10\} \to {\mathbb Z}_+,
$$
setting $\xi(M) = \frac{(M-9)(M-8)}{2} + 12$.

For simplicity of notations, we identify a pair of irreducible
polynomials $(g,h)\in {\cal F}$ with the corresponding Fano double
cover $V$ and write $V\in {\cal F}$; this can not lead to any
confusion. Now we can state the main result of the
paper.\vspace{0.1cm}

{\bf Theorem 1.} {\it There exists a Zariski open subset ${\cal
F}_{\rm reg}\subset {\cal F}$ such that the following claims are
true.\vspace{0.1cm}

{\rm (i)} Every variety $V\in {\cal F}_{\rm reg}$ is factorial and
has at most terminal singularities.\vspace{0.1cm}

{\rm (ii)} The complement ${\cal F}\setminus {\cal F}_{\rm reg}$
is of codimension at least $\xi(M)$ in ${\cal F}$.\vspace{0.1cm}

{\rm (iii)} Every variety $V\in {\cal F}_{\rm reg}$ is
birationally superrigid.}\vspace{0.1cm}

{\bf Corollary 1.} {\it For every variety $V\in {\cal F}_{\rm
reg}$ the following claims are true.\vspace{0.1cm}

{\rm (i)} Every birational map $V\dashrightarrow V'$ to a Fano
variety with ${\mathbb Q}$-factorial terminal singularities and
Picard number 1 is a biregular isomorphism.\vspace{0.1cm}

{\rm (ii)} There are no rational dominant maps $V\dashrightarrow
S$ onto a positive-dimensional variety $S$, the general fibre of
which is rationally connected (or has negative Kodaira dimension).
In particular, there are no structures of a Mori fibre space over
a positive-dimensional base on $V$.\vspace{0.1cm}

{\rm (iii)} The variety $V$ is non-rational and its groups of
birational and biregular automorphisms are the same:} $\mathop{\rm
Bir} V=\mathop{\rm Aut} V$.\vspace{0.1cm}

{\bf Proof of the corollary.} The claims (i)-(iii) are all the
standard implications of the property of being birationally
superrigid, see, for instance, \cite[Chapter 2]{Pukh13a}.
Q.E.D.\vspace{0.3cm}


{\bf 0.2. The regularity conditions.} The open subset ${\cal
F}_{\rm reg}\subset {\cal F}$ is defined by a number of explicit
local conditions, to be satisfied at every point, which we now
list. Let $o\in V$ be a point, $p=\sigma(o)\in G$ its image on
${\mathbb P}$. We assume, therefore, that $g(p)=0$. Let
$z_1,\dots, z_{M+1}$ be a system of affine coordinates on
${\mathbb P}$ with the origin at $p$ and
$$
g=q_1+q_2+\dots +q_m,\quad h=w_0+w_1+w_2+\dots +w_{2l}
$$
the decomposition of $g,h$ (dehomogenized but for simplicity of
notations denoted by the same symbols) into components,
homogeneous in $z_*$. We may assume that $z_i=x_i/x_0$ are
coordinates on the affine chart $\{x_0\neq 0\}$ on ${\mathbb P}$.
Adding the new affine coordinate $y=u/x_0^l$, we extend that chart
to
$$
{\mathbb A}^{M+2}_{z_*,y}\subset {\mathbb P}(1^{M+2},l),
$$
where the variety $V$ is a complete intersection, given by the
system of two equations:
$$
\begin{array}{c}
\phantom{-y^2+w_0+}q_1+q_2+\dots +q_m=0,\\
-y^2+w_0+w_1+w_2+\dots +w_{2l}=0.
\end{array}
$$
Note that $p\in W$ if and only if $w_0=0$.\vspace{0.1cm}

We assume that the hypersurface $G\subset {\mathbb P}$ has at most
quadratic singularities: if $q_1\equiv 0$, then $q_2\not\equiv 0$.
Furthermore, we assume that $G$ is {\it regular} in the standard
sense at very point $p\in G$:\vspace{0.1cm}

(R0.1) If $q_1\not\equiv 0$, then the sequence
$$
q_1,\, q_2, \dots, q_m
$$
is regular in ${\cal O}_{p,{\mathbb P}}$.\vspace{0.1cm}

(R0.2) If $q_1\equiv 0$, then the sequence
$$
\phantom{q_1,\,} q_2, \dots, q_m
$$
is regular in ${\cal O}_{p,{\mathbb P}}$.\vspace{0.1cm}

We will need also some additional regularity conditions for the
polynomials $g,h$ at the point $p$, which depend on whether $p\in
W$ or $p\not\in W$ and on the type of singularity $o\in V$ that we
allow.\vspace{0.1cm}

We start with the {\bf non-singular case}.\vspace{0.1cm}

(R1.1) If $w_0\neq 0$, then we have no additional conditions (only
(R0.1) is needed).\vspace{0.1cm}

(R1.2) If $w_0=0$, then
$$
q_2|_{\{q_1=w_1=0\}}\not\equiv 0.
$$
Note that in the second case as the point $o\in V$ is assumed to
be non-singular, the linear forms $q_1$ and $w_1$ must be linearly
independent.\vspace{0.1cm}

Now let us consider the {\bf quadratic case}.\vspace{0.1cm}

Here we have three possible ways of getting a singular point and,
accordingly, three types of regularity conditions.\vspace{0.1cm}

(R2.1) Out side the ramification divisor: if $w_0\neq 0$, then
$q_1\equiv 0$ and
$$
\mathop{\rm rk} q_2\geqslant 7.
$$

(R2.2) On the ramification divisor with $G$ non-singular: $w_0=0$,
$q_1\not\equiv 0$, $w_1\equiv 0$ and
$$
\mathop{\rm rk} w_2|_{\{q_1=0\}}\geqslant 6.
$$

(R2.3) On the ramification divisor with $G$ singular: $w_0=0$,
$q_1\equiv 0$, $w_1\not\equiv 0$ and
$$
\mathop{\rm rk} q_2|_{\{w_1=0\}}\geqslant 7.
$$

Apart from non-singular points and quadratic singularities, we
allow more complicated points which we call {\bf bi-quadratic}.
Assume that $w_0=0$ and $q_1\equiv w_1\equiv 0$.\vspace{0.1cm}

(R$2^2$) For a general $11-$dimensional linear subspace $P\subset
{\mathbb C}^{M+2}_{z_*,y}$ the closed algebraic set
$$
Q_P=\left\{q_2|_P=(y^2-w_2)|_P=0\right\}\subset {\mathbb
P}(P)\cong {\mathbb P}^{10}
$$
is a non-singular complete intersection of codimension
2.\vspace{0.1cm}

We say that a pair $(g,h)\in {\cal F}$ is {\it regular} if the
hypersurface $G=\{g=0\}\subset {\mathbb P}$ is regular at every
point in the sense of the conditions (R0.1) and (R0.2) (whichever
applies at the given point), and the relevant regularity condition
from the list above is satisfied at every point $o\in
\{g=(u^2-h=0\}$.\vspace{0.1cm}

Note that $(g,h)\in {\cal F}$ being regular implies that the
closed set
$$
V=\{g=u^2-h=0\}\subset {\mathbb P}(1^{M+2},l)
$$
is an irreducible complete intersection of codimension 2, the
singular points of which are either quadratic singularities of
rank $\geqslant 7$ or bi-quadratic singularities satisfying the
condition (R$2^2$). In any case, the singularities of $V$ are
complete intersection singularities and the singular locus
$\mathop{\rm Sing} V$ has codimension at least 7 in $V$, so the
Grothendieck theorem on parafactoriality \cite{CL} applies and $V$
turns out to be a factorial variety. Furthermore, it is easy to
check that the property of having at most quadratic singularities
of rank $\geqslant r$ is stable with respect to blowing up
non-singular subvarieties (see \cite[Section 3.1]{Pukh15d} for a
detailed proof and discussion, and the same arguments apply to
bi-quadratic singularities satisfying (R$2^2$)), so that, in
particular, the singularities of $V$ are terminal.\vspace{0.1cm}

Now setting ${\cal F}_{\rm reg}\subset {\cal F}$ to be the open
subset of regular pairs $(g,h)$ (or, abusing the notations,
regular varieties $V=V(g,h)$), we get the claim (i) of Theorem
1.\vspace{0.1cm}

Therefore, Theorem 1 is implied by the following two
claims.\vspace{0.1cm}

{\bf Theorem 2.} {\it The complement ${\cal F}\setminus {\cal
F}_{\rm reg}$ is of codimension at least $\xi(M)$ in} ${\cal
F}$.\vspace{0.1cm}

{\bf Theorem 3.} {\it A regular variety $V\in {\cal F}_{\rm reg}$
is birationally superrigid.}\vspace{0.3cm}


{\bf 0.3. The structure of the paper.} We prove Theorem 3 in
Section 1 and Theorem 2 in Section 2. The arguments are
independent of each other.\vspace{0.1cm}

In order to prove Theorem 3, we assume the converse: $V$ is {\it
not} birationally superrigid. This implies, in a standard way
\cite[Chapter 2, Section 1]{Pukh13a} that there is a mobile linear
system $\Sigma\subset |nH|$ with a {\it maximal singularity}. The
centre of the maximal singularity is an irreducible subvariety
$B\subset V$.There are a number of options for $B$: it can have a
small ($\leqslant 4$) codimension or a higher ($\geqslant 5$)
codimension in $V$, be contained or not contained in the singular
locus $\mathop{\rm Sing} V$ (and more specifically, in the locus
of bi-quadratic points), be contained or not contained in the
ramification divisor. For each of these options, we {\it exclude}
the maximal singularity, that is, we show that its existence leads
to a contradiction. After that, we conclude that the initial
assumption was incorrect and $V$ is birationally
superrigid.\vspace{0.1cm}

Theorem 2 is shown by different and very explicit arguments. We
fix a point $o\in {\mathbb P}(1^{M+2},l)$ and consider varieties $V\ni o$.
For each type of the point $o$ (from the list given in Subsection
0.2) and each regularity condition we estimate the codimension of
the closed set of pairs $(g,h)\in {\cal F}$ such that
$o\in\{g=u^2-h=0\}$ and the condition under consideration is
violated. Taking the minimum of our estimates, we prove Theorem
2.\vspace{0.1cm}

The decisive point of this paper is applying the generalized
$4n^2$-{\it inequality} \cite{Pukh2017b} to excluding the maximal
singularities, the centre of which is contained in the quadratic
or bi-quadratic locus: without it, the task would have been too
hard. The regularity conditions make sure that the generalized
$4n^2$-inequality applies. Given the new essential ingredient,
excluding the maximal singularity becomes
straightforward.\vspace{0.3cm}


{\bf 0.4. Historical remarks.} We say that a theorem stating
birational (super)rigidity is {\it effective}, if it contains an
effective bound for the codimension of the set of non-rigid
varieties (in the natural parameter space of the family under
consideration). The first effective result was obtained in
\cite{EP1}. For complete intersections see \cite{EvP2017,EvP2018}.
The importance of effective results is explained by the problem of
birational rigidity of Fano-Mori fibre spaces, see \cite{Pukh15d},
generalizing the famous Sarkisov theorem \cite{S82} to fibre
spaces with higher-dimensional fibres.\vspace{0.1cm}

Birational rigidity of certain mildly singular Fano double covers
was shown in \cite{Ch06a,Johnstone}. The result of
\cite{Johnstone} was effective in our sense. Iterated double
covers and cyclic covers of degree $\geqslant 3$ were considered
in \cite{Pukh03} and \cite{Pukh09a}, respectively (only
non-singular varieties were treated in these papers). Triple
covers with singularities were shown to be birationally superrigid
in \cite{Ch04c}. For a study of the question, how many families of
higher-dimensional non-singular Fano complete intersections are
there in the weighted complete intersections, see
\cite{PrzShr2016}.\vspace{0.3cm}


{\bf 0.5. Acknowledgements.} The second author is grateful to the
Leverhulme Trust for the financial support (Research Project Grant
RPG-2016-279).


\section{Proof of birational superrigidity}

In this section we prove Theorem 3. First, we remind the
definition and some basic facts about maximal singularities,
classifying them and excluding the cases of low codimension of the
centre (Subsection 1.1). Then we exclude the maximal
singularities, the centre of which is not contained in the
singular locus of $V$ (Subsection 1.2). Finally, we exclude the
cases when the centre of a maximal singularity is contained in the
singular locus (Subsection 1.3). The last group of cases, which
traditionally was among the hardest to deal with, now becomes the
easiest due to the generalized $4n^2$-inequality shown in
\cite{Pukh2017b}.\vspace{0.3cm}

{\bf 1.1. Maximal singularities.}Assume that a fixed regular
double hypersurface $V\in {\cal F}_{\rm reg}$ is not birationally
superrigid. It is well known (see, for instance, \cite[Chapter 2,
Section 1]{Pukh13a}), that this assumption implies that there is a
mobile linear system $\Sigma\subset |nH|$, a birational morphism
$\varphi\colon\widetilde{V} \to V$ and a $\varphi$-exceptional
prime divisor $E\subset \widetilde{V}$, satisfying the {\it
Noether-Fano inequality}
$$
\mathop{\rm ord}\nolimits_E \varphi^*\Sigma > n\cdot a(E).
$$
Here $\widetilde{V}$ is assumed to be non-singular projective,
$\varphi$ a composition of blow ups with non-singular centres,
$a(E)=a(E,V)$ is the discrepancy of $E$ with respect to $V$. The
prime divisor $E$ (or the discrete valuation of the field of
rational functions ${\mathbb C}(\widetilde{V})\cong {\mathbb
C}(V)$) is called a {\it maximal singularity} of the system
$\Sigma$. Equivalently, for any divisor $D\in Sigma$ the pair
$(V,\frac{1}{n} D)$ is not canonical with $E$ a non-canonical
singularity of the pair. Set $B=\varphi(E)\subset V$ to be its
centre on $V$ and $\overline{B}=\sigma(B)\subset G$ its projection
on ${\mathbb P}$. We have the following options:\vspace{0.1cm}

(1) $\mathop{\rm codim} (B\subset V)=2$,\vspace{0.1cm}

(2) $\mathop{\rm codim} (B\subset V)=3$ or $4$,\vspace{0.1cm}

(3) $\mathop{\rm codim} (B\subset V) \geqslant 5$ and
$\overline{B}\not\subset W$, $B\not\subset \mathop{\rm Sing}
V$,\vspace{0.1cm}

(4) $\mathop{\rm codim} (B\subset V)\geqslant 5$ and
$\overline{B}\subset W$, $B\not\subset \mathop{\rm Sing}
V$,\vspace{0.1cm}

(5) $B$ is contained in the (closure of the) locus of quadratic
singularities, but not in the locus of bi-quadratic
singularities,\vspace{0.1cm}

(6) $B$ is contained in the locus of bi-quadratic
singularities.\vspace{0.1cm}

We have to show that none of these cases take place. Note that the
inequality
\begin{equation}\label{27.09.2018.1}
\mathop{\rm mult}\nolimits_B \Sigma >n
\end{equation}
holds. Let $Z=(D_1\circ D_2)$ be the algebraic cycle of
scheme-theoretic intersection of general divisors $D_1,D_2\in
\Sigma$, the {\it self-intersection} of the system $\Sigma$. Note
that $Z\sim n^2 H^2$.\vspace{0.1cm}

Our first observation is that the case (1) does not realize.
Indeed, let $P$ be a general $7$-dimensional plane in ${\mathbb
P}$. Then $V_P=V\cap \sigma^{-1}(P)$ is a non-singular
$6$-dimensional variety. By the Lefschetz theorem,
$$
\mathop{\rm Pic} V_P ={\mathbb Z} H_P\quad\mbox{and}\quad A^2 V_P=
{\mathbb Z} H^2_P,
$$
where $H_P$ is the hyperplane section and $A^2$ the numerical Chow
group of codimension 2 cycles. The restriction $Z_P=(Z\circ
V_P)\sim n^2 H_P$ is an effective cycle. If $\mathop{\rm codim}
(B\subset V)=2$, then $Z$ contains $B$ as a component with
multiplicity at least $(\mathop{\rm mult}\nolimits_B \Sigma)^2$;
therefore, $Z_P$ contains $B_P=(B\circ V_P)=B\cap V_P$ with
multiplicity at least $(\mathop{\rm mult}\nolimits_B \Sigma)^2$.
However, $B_P\sim b H^2_P$ for some $b\geqslant 1$ and the
inequality (\ref{27.09.2018.1}) can not be true. So we may assume
that $\mathop{\rm codim} (B\subset V)\geqslant 3$.\vspace{0.1cm}

{\bf Proposition 1.1.} {\it The case (2) does not
realize.}\vspace{0.1cm}

{\bf Proof.} Assume the converse: $\mathop{\rm codim} (B\subset
V)\in\{3,4\}$. Then $B\not\subset \mathop{\rm Sing} V$ and so the
standard $4n^2$-inequality holds:
$$
\mathop{\rm mult}\nolimits_B Z>4n^2,
$$
see \cite[Chapter 2]{Pukh13a}. Again, take a general
$7$-dimensional plane $P\subset {\mathbb P}$ and let $V_P$, $Z_P$,
$H_P$ and $B_P$ mean the same as above. We can find an irreducible
subvariety $Y\sim d H^2_P$ of codimension 2 in $V_P$ such that
$$
\mathop{\rm mult}\nolimits_{B_P} Y> 4d.
$$
Set $G_P=G\cap P$: it is a non-singular hypersurface of degree $m$
in $P\cong {\mathbb P}^7$. Writing $H_G$ for the class of its
hyperplane section, we get
$$
\mathop{\rm Pic} G_P ={\mathbb Z} H_G\quad\mbox{and}\quad A^2 G_P=
{\mathbb Z} H^2_G.
$$
Let $\overline{Y}=\sigma(Y)\subset G_P$ and
$\overline{B}_P=\sigma(B_P)$ be the images of $Y$ and $B_P$,
respectively. Then
$$
\overline{Y}\sim d^* H^2_G
$$
with $d^*=d$ or $\frac12 d$, and the inequality
$$
\mathop{\rm mult}\nolimits_{\overline{B}_P}
\overline{Y} > 2d^*
$$
holds. But $\mathop{\rm dim} \overline{B}_P\in \{2,3\}$, so we get
a contradiction with \cite[Proposition 5]{Pukh02f} (see also
``Pukhlikov's Lemma'' in \cite{Suzuki}). Q.E.D. for Proposition
1.1.\vspace{0.1cm}

From now on, we assume that $\mathop{\rm codim} (B\subset
V)\geqslant 5$.\vspace{0.1cm}

In order to exclude the cases (3-6), we will need the regularity
conditions (R0.1,2), or rather, the facts that are summarized in
the proposition below.\vspace{0.1cm}

{\bf Proposition 1.2} {\it Let $S\subset G$ be an irreducible
subvariety of codimension $a\in \{2,3\}$ and $p\in S$ a
point.\vspace{0.1cm}

{\rm (i)} Assume that $G$ is non-singular at $p$. Then
$$
\mathop{\rm mult}\nolimits_p S\leqslant \frac{a+1}{m} \mathop{\rm
deg} S.
$$

{\rm (ii)} Assume that $G$ is singular at $p$. Then}
$$
\mathop{\rm mult}\nolimits_p S\leqslant \frac{a+2}{m} \mathop{\rm
deg} S.
$$

{\bf Proof.} The claims are the standard implications of the
regularity conditions (R0.1,2). see, for instance, \cite[Chapter
3]{Pukh13a} for the standard arguments delivering the estimates
for the multiplicity in terms of degree. Q.E.D.\vspace{0.3cm}


{\bf 1.2. The non-singular case.} Let us exclude the options (3)
and (4). Here $\overline{B}\not\in \mathop{\rm Sing} G$ and in any
case $\overline{B}\not\in \mathop{\rm Sing} W$.\vspace{0.1cm}

{\bf Proposition 1.3.} {\it The case (4) does not
realize.}\vspace{0.1cm}

{\bf Proof.} Here we can argue in word for word the same way as in
\cite[Subsection 3.3, Case 2]{Pukh00e}: take a general point $o\in
B$, so that $p=\sigma(o)\in W$ is a non-singular point on $W$. The
tangent hyperplanes
$$
T_p G\quad \mbox{and}\quad T_p W_{\mathbb P}
$$
are distinct and their $\sigma$-preimages on $V$ are singular.
Therefore,
$$
\Delta=\sigma^{-1}\left(T_pG\cap T_p W_{\mathbb P} \cap G \right)
$$
is an irreducible subvariety of codimension 2 on $V$, satisfying
the relations
$$
\Delta\sim H^2\quad\mbox{and}\quad\mathop{\rm mult}\nolimits_o
\Delta=4,
$$
the second equality is guaranteed by the regularity condition
(R1.2).\vspace{0.1cm}

On the other hand, from the (standard) $4n^2$-inequality we get
that there is an irreducible subvariety $Y\subset V$ such that
$$
Y\sim dH^2\quad\mbox{and}\quad \mathop{\rm mult}\nolimits_o Y > 4d
$$
for some $d\in {\mathbb Z}_+$. Therefore, $Y\neq\Delta$, which
means that $Y$ is not contained in at least one of the two
divisors
$$
\sigma^{-1}(T_p G)\quad \mbox{and}\quad \sigma^{-1}(T_p W_{\mathbb
P}).
$$
Taking the scheme-theoretic intersection of $Y$ with that divisor
and selecting a suitable irreducible component, we obtain an
irreducible subvariety $Y^*\subset V$ of codimension 3 such that
$$
\mathop{\rm mult}\nolimits_o Y^* >\frac{4}{m} \mathop{\rm
deg}\nolimits_H Y^*.
$$
The image $S=\sigma(Y^*)\subset G$ is an irreducible subvariety of
codimension 3, satisfying the inequality
$$
\mathop{\rm mult}\nolimits_p S > \frac{4}{m} \mathop{\rm deg} S.
$$
We get a contradiction with the claim (i) of Proposition 1.2.
Q.E.D.\vspace{0.1cm}

{\bf Proposition 1.4.} {\it The case (3) does not
realize.}\vspace{0.1cm}

{\bf Proof.} Assume the converse. Let $o\in B$ be a general point,
so that $p=\sigma(o)\not\in W$ and $p\not\in \mathop{\rm Sing} V$.
Note that $\sigma_*\colon T_oV\to T_pG$ is an isomorphism of
vector spaces. Let $\lambda\colon V^+\to V$ be the blow up of the
point $o$ and $\lambda_G\colon G^+\to G$ the blow up of the point
$p$, with the exceptional divisors $E^+$ and $E^+_G$,
respectively. We have the natural isomorphism
$$
E^+\stackrel{\sigma}{\longrightarrow} E^+_G\cong {\mathbb
P}^{M-1}.
$$
It is well known (the ``$8n^2$-inequality'', see, for instance,
\cite[Chapter 2]{Pukh13a}), that there is a linear subspace
$\Lambda\subset E^+$ of codimension 2 such that
$$
\mathop{\rm mult}\nolimits_o Z + \mathop{\rm
mult}\nolimits_{\Lambda} Z^+ > 8n^2,
$$
where $Z^+$ is the strict transform of the self-intersection $Z$
on $V^+$. Let $P\subset {\mathbb P}$ be a general hyperplane such
that
$$
\sigma^{-1}(G\cap P)^+\supset\Lambda.
$$
Set $G_P=G\cap P$; obviously, for a general $P$ none of the
irreducible components of $Z$ is contained in $\sigma^{-1}(G_P)$.
Therefore, the effective cycle
$$
Z_P=(Z\circ G_P)
$$
of codimension 3 on $V$ satisfies the inequality
$$
\mathop{\rm mult}\nolimits_o Z_P > 8n^2.
$$
Taking a suitable irreducible component $Y$ of $Z_P$ and its image
$S=\sigma(Y)$, we obtain an irreducible subvariety $S\subset G$ of
codimension 3, satisfying at the non-singular point $p\in G$ the
inequality
$$
\mathop{\rm mult}\nolimits_p S > \frac{4}{m} \mathop{\rm deg} S.
$$
This contradicts the claim (i) of Proposition 1.2.
Q.E.D.\vspace{0.3cm}


{\bf 1.3. The singular case.} It remains to exclude the options
(5) and (6), where $B\subset \mathop{\rm Sing} V$. It is here that
we use the generalized $4n^2$-inequality shown in
\cite{Pukh2017b}.\vspace{0.1cm}

{\bf Proposition 1.5.} {\it The cases (5) and (6) do not
realize.}\vspace{0.1cm}

{\bf Proof.} Assume that the case (5) takes place. Let $o\in B$ be
a point of general position. The singularity $o\in V$ is a
quadratic singularity, satisfying the requirements of the main
theorem of \cite{Pukh2017b}. Therefore,
$$
\mathop{\rm mult}\nolimits_o Z > 4n^2\cdot \mathop{\rm
mult}\nolimits_o V=8n^2.
$$
Taking a suitable irreducible component $Y$ of $Z$ and its image
$S=\sigma(Y)$, we obtain an irreducible subvariety $S\subset G$ of
codimension 2, satisfying at the quadratic point $p\in G$ the
inequality
$$
\mathop{\rm mult}\nolimits_p S > \frac{4}{m} \mathop{\rm deg} S,
$$
which contradicts the claim (ii) of Proposition 1.2.\vspace{0.1cm}

The case (6) is excluded in a similar way, just for $Z$ we get the
inequality
$$
\mathop{\rm mult}\nolimits_o Z > 4n^2\cdot \mathop{\rm
mult}\nolimits_o V=16n^2.
$$
and for $S$ the inequality
$$
\mathop{\rm mult}\nolimits_p S > \frac{8}{m} \mathop{\rm deg} S,
$$
which can not be satisfied at a quadratic point $p\in G$ by
Proposition 1.2. Q.E.D.\vspace{0.1cm}

Proof of Theorem 3 is now complete.


\section{Estimates for the codimension}

In this section we prove Theorem 3. To this purpose, for each $M \geqslant 10$ we construct an algebraic subset ${\cal Z} \subset {\cal F}$ of codimension $\geqslant \xi(M)$, such that ${\cal F} - {\cal Z} \subset {\cal F}_{\rm reg}$.

As a first step we reduce the construction to double hypersurfaces containing a fixed point $o\in {\mathbb P}(1^{M+2},l)$: The point $[(0: \cdots : 0) :_l 1] \in {\mathbb P}(1^{M+2},l)$ is contained in no such double hypersurface, by its construction. For all other points $o = [o^\prime :_l u] \in {\mathbb P}(1^{M+2},l)$ the subset ${\cal F}^o \subset {\cal F}$ of pairs $(g,h) \in {\cal{F}}$ such that $o$ is contained in
$$
V = \{g=0, u^2 = h\} \subset {\mathbb P}(1^{M+2},l),
$$
the double cover of $G = \{g = 0\}$ associated to $(g,h)$, is equal to ${\cal P}_{m,M+2}^o \times {\cal P}_{2l,M+2}^o$, with
$$
{\cal P}_{m,M+2}^o = \{g \in {\cal P}_{m,M+2}: g(o^\prime)=0\}\ \mathrm{and\ }
{\cal P}_{2l,M+2}^o = \{h \in {\cal P}_{2l,M+2}: u^2 = h(o^\prime)\}
$$
affine hyperplanes of ${\cal P}_{m,M+2}$ resp.\ ${\cal P}_{2l,M+2}$.

Now choose a point $o_1 = [o_1^\prime :_l u_1] \in {\mathbb P}(1^{M+2},l) \setminus \{[0: \cdots : 0 :_l 1]\}$ with $u_1 \neq 0$ and a point $o_2 = [o_2^\prime :_l 0] \in {\mathbb P}(1^{M+2},l)$. \vspace{0.1cm}

{\bf Proposition 2.1.} {\it For $i = 1, 2$ let ${\cal Z}_{o_i} \subset {\cal F}^{o_i}$ be algebraic subsets such that ${\cal F}^{o_i} \setminus {\cal Z}_i \subset {\cal F}^{o_i}_{\rm reg}$. Then there exists an algebraic subset ${\cal Z} \subset {\cal F}$ such that ${\cal F} \setminus {\cal Z} \subset {\cal F}_{\rm reg}$ and}
$$
{\rm codim}_{\cal F} {\cal Z} \geqslant \min({\rm codim}_{{\cal F}^{o_1}} {\cal Z}_{o_1} - M, {\rm codim}_{{\cal F}^{o_2}} {\cal Z}_{o_2} - M + 1).
$$

{\bf Proof.} ${\rm PGL}(M+2)$ acts on ${\mathbb P}(1^{M+2},l)$ by transforming the first $M+2$ homogeneous coordinates in the standard way. This action has the three orbits $\{[0: \cdots : 0 :_l 1]\}$, $\{[o^\prime :_l u] \in {\mathbb P}(1^{M+2},l)  : u \neq 0\} \setminus \{[0: \cdots : 0 :_l 1]\}$ and $\{[o^\prime :_l u] \in {\mathbb P}(1^{M+2},l): u = 0\}$. Thus, for each point $o_1 \in \{u \neq 0\} \setminus \{[0: \cdots : 0 :_l 1]\}$ resp.\ $o_2 \in \{u = 0\}$ we can find isomorphic algebraic subset ${\cal Z}_{o_1}$ resp.\ ${\cal Z}_{o_2}$ such that  ${\cal F}^{o_1} \setminus {\cal Z}_{o_1} \subset {\cal F}^{o_1}_{\rm reg}$ resp.\ ${\cal F}^{o_2} \setminus {\cal Z}_{o_2} \subset {\cal F}^{o_2}_{\rm reg}$.

The closure ${\cal Z}_1$ of the union of all the ${\cal Z}_{o_1}$ has dimension $\leqslant \dim {\cal Z}_{o_1} + M + 2$, whereas the closure ${\cal Z}_2$ of the union of all the ${\cal Z}_{o_2}$ has dimension $\leqslant \dim {\cal Z}_2 + M + 1$. Since ${\rm codim}_{{\cal F}} {\cal F}^o = 2$ this implies the bound on the codimension of ${\cal Z} = {\cal Z}_1 \cup {\cal Z}_2$. Q.E.D.\vspace{0.1cm}

Note that a point $o \in \{u \neq 0\}$ can only lie outside the ramification locus of a Fano double cover $V$, whereas a point $o \in \{u=0\}$ must lie on the ramification locus. \vspace{0.3cm}


{\bf 2.1. Codimension estimates for points outside the ramification locus.} Choose a point $o \in \{u \neq 0\} \setminus \{[0: \cdots : 0 :_l 1]\}$. We first treat the cases when the regularity conditions on the hypersurface $G=\{g=0\}\subset {\mathbb P}$ fail.

Using the notation in the Introduction assume that $q_1 \not\equiv 0$. The set $S_{R0.1}$ of pairs $(g,h)$ in ${\cal F}^o$ such that $q_1, \ldots, q_m$ is not a regular sequence in ${\cal O}_{p,{\mathbb P}}$ is a closed algebraic subset of the Zariski-open subset $\{q_1 \not\equiv 0\} \subset {\cal F}^o$. It is stratified according to the position where $q_1, \ldots, q_m$ is not any longer regular: Since $q_1 \not\equiv 0$ this can only happen from $q_2$ on, so $S_{R0.1} = S_{R0.1}^2 \cup \ldots \cup S_{R0.1}^m$ with
$$
S_{R0.1}^d = \{(q_1, \ldots, q_m; h): q_2, \ldots, q_{d-1} {\rm \ is\ regular,\ but\ not\ } q_2, \ldots, q_d \} \subset {\cal F}^o.
$$
for $d = 2, \ldots, m$. The set $S_{R0.1}^d$ is closed algebraic in $S_{R0.1} \setminus \bigcup_{i=2}^{d-1} S_{R0.1}^i$, thus the codimension of its Zariski closure in ${\cal F}^o$ is $\geqslant$ to the codimension of its intersection with the fiber in ${\cal F}^o$ over a fixed regular sequence $q_1, \ldots, q_{d-1}$ in this fiber, under the natural projection. By the methods in \cite{Pukh98b} this codimension is $\geqslant \binom{M+1}{d}$ for $2 \leq d \leq m$. Since $m+2l \leqslant M+1$ this implies:
\begin{equation}\label{codimA1-eq}
{\rm codim}_{{\cal F}^o} S_{R0.1} \geqslant \binom{M+1}{2}.
\end{equation}

If $q_1 \equiv 0$ and $S_{R0.2}$ denotes the set of pairs $(g,h)$ in ${\cal F}^o$ such that $q_2, \ldots, q_m$ is not a regular sequence in ${\cal O}_{p,{\mathbb P}}$, we find as before a lower bound for the codimension of the closed algebraic subset $S_{R0.2}$ in ${\cal F}^o$:
\begin{equation}\label{codimA2-eq}
{\rm codim}_{{\cal F}^o} S_{R0.2} \geqslant M+1 + \binom{M+1}{2} = \binom{M+2}{2}.
\end{equation}
Here, the summand $M+1$ counts the codimensions given by the vanishing of $q_1$.

Next, we study the case when the point $o$ is too singular on the double cover $V$, that is when condition (R2.1) fails. This happens when $q_1 \equiv 0$ and $\mathop{\rm rk} q_2 \leqslant 6$, and we denote the closed algebraic subset of pairs $(g,h)$ in ${\cal F}^o$ satisfying these conditions by $S_{R2.1}$.

Quadratic forms in $M+1$ variables correspond to symmetric $(M+1) \times (M+1)$ matrices parametrised by a $\binom{M+2}{2}$-dimensional affine space ${\rm Sym}_{M+1}$, and the rank of a quadratic form $q_2$ equals the rank of the corresponding symmetric matrix $A$. But $\mathop{\rm rk} A \leqslant r$ if and only if there exists an $(M+1-r)$-dimensional vector subspace $\Lambda \subset \mathbb{C}^{M+1}$ spanned by $0$-eigenvectors of $A$. Such matrices $A \in {\rm Sym}_{M+1, \leqslant r} = \{A \in {\rm Sym}_{M+1} : \mathop{\rm rk} A \leqslant r\}$ lie in the image of the incidence variety
$$
\Phi = \{(A, \Lambda) : A \cdot v = 0\ {\rm for\ all\ } v \in \Lambda\} \subset {\rm Sym}_{M+1} \times {\rm Gr}(M+1-r, M+1)
$$
under the projection to ${\rm Sym}_{M+1}$. This projection has $0$-dimensional general fibers, for matrices of rank $r$, so
${\rm codim}_{{\rm Sym}_{M+1}} {\rm Sym}_{M+1, \leqslant r} = \dim {\rm Sym}_{M+1} - \dim \Phi$. On the other hand, the projection of $\Phi$ onto the Grassmannian ${\rm Gr}(M+1-r, M+1)$ has fibers of dimension $\binom{r+1}{2}$, so $\dim \Phi = \binom{r+1}{2} + r(M+1-r)$ and
$$
{\rm codim}_{{\rm Sym}_{M+1}} {\rm Sym}_{M+1, \leqslant r} = \binom{M+2}{2} - \binom{r+1}{2} - r(M+1-r) =
\frac{(M+2-r)(M+1-r)}{2}.
$$

Setting $r=6$ and adding the $M+1$ codimensions given by $q_1 \equiv 0$ we obtain
\begin{equation}\label{codimA3-eq}
{\rm codim}_{{\cal F}^o} S_{R2.1} =  M+1 + \frac{(M-4)(M-5)}{2} = \frac{(M-4)(M-3)}{2} + 5.
\end{equation}


{\bf 2.2. Codimension estimates for points on the ramification locus.} Choose a point $o \in \{u = 0\}$. Using the notation in the Introduction $o$ will lie on a double cover given by a pair $(g,h) \in {\cal F}^o$ only if $w_0 = 0$.

As for points outside the ramification locus we obtain the following two codimension bounds for subsets $S_{R0.1} \subset {\cal F}^o$ and $S_{R0.2} \subset {\cal F}^o$ where the regularity conditions on the hypersurface $G=\{g=0\}\subset {\mathbb P}$ fail:
\begin{equation}\label{codimB1-eq}
{\rm codim}_{{\cal F}^o} S_{R0.1} \geqslant \binom{M+1}{2}
\end{equation}
and
\begin{equation}\label{codimB2-eq}
{\rm codim}_{{\cal F}^o} S_{R0.2} \geqslant \binom{M+2}{2}.
\end{equation}

Next, we study the set $S_{\rm R1.2} \subset {\cal F}^o$ of pairs $(g,h)$ such that $o$ is non-singular on the associated double cover but condition (R1.2) fails. That is the case when $q_1 \not\equiv 0$, $w_1 \not\equiv \lambda q_1$ for all $\lambda \in \mathbb{C}$ and $q_{2|\{q_1=w_1=0\}} \equiv 0$. The last identity is equivalent to $q_2 \equiv q_1 \cdot q_1^\prime + w_1 \cdot w_1^\prime$ for two linear forms $q_1^\prime, w_1^\prime$. Since the first two conditions are open in ${\cal F}^o$ it is enough to determine the codimension of the set of $q_2$ in the space of all quadratic forms in $M+1$ variables that are of the above form for given $q_1, w_1$: By a change of coordinates $q_1$ and $w_1$ may be identified with two of the $M+1$ variables, thus the requested codimension equals the dimension of quadratic forms in $M-1$ variables. So we have
\begin{equation}\label{codimB2-eq}
{\rm codim}_{{\cal F}^o} S_{R1.2} \geqslant \binom{M}{2}.
\end{equation}

Pairs $(g,h)$ for which $o$ is a singular point on the associated double cover mapped to a non-singular point on the hypersurface $G \subset \mathbb{P}$ fail condition (R2.2) if and only if $q_1 \not\equiv 0$, $w_1 = \lambda q_1$ for some $\lambda \in \mathbb{C}$ and ${\rm rk}(w_2 - \lambda q_{2|\{q_1 = 0\}}) \leqslant 5$. The codimension of the set $S_{\rm R2.2} \subset {\cal F}^o$ of such pairs equals the sum of $M$ (from $w_1 \equiv \lambda q_1$) and the codimension of quadratic forms of rank $\leqslant 5$ when restricted to a given linear form, in the space of all quadratic forms in $M+1$ variables. Since by a coordinate change we can assume that $q_1$ is one of the $M+1$ variables it is enough to calculate the codimension of quadratic forms of rank $\leqslant 5$ in the space of all quadratic forms in $M$ variables. Imitating the calculations in Section 2.1 we obtain a lower bound for this codimension as
$\frac{(M-4)(M-5)}{2}$.
Adding up this leads to
\begin{equation}\label{codimB2-eq}
{\rm codim}_{{\cal F}^o} S_{R2.2} \geqslant M + \frac{(M-4)(M-5)}{2} = \frac{(M-4)(M-3)}{2} + 4.
\end{equation}

Pairs $(g,h)$ for which $o$ is a singular point on the associated double cover mapped to a singular point on the hypersurface $G \subset \mathbb{P}$ fail condition (R2.3) if and only if $q_1 \equiv 0$, $w_1 \not\equiv 0$ and ${\rm rk}(q_{2|\{w_1 = 0\}}) \leqslant 6$. As before we obtain a lower bound for the codimension of the set $S_{\rm R2.3} \subset {\cal F}^o$ of such pairs as
\begin{equation}\label{codimB2-eq}
{\rm codim}_{{\cal F}^o} S_{R2.3} \geqslant (M + 1) + \frac{(M-5)(M-6)}{2} = \frac{(M-5)(M-4)}{2} + 6.
\end{equation}

Finally, we need to look at the set $S_{{\rm R2.2}^{\rm 2}} \subset {\cal F}^o$ of pairs $(g,h)$ where $o$ is a biquadratic singular point on the associated double cover failing condition ${\rm R2.2}^{\rm 2}$. This is the case if and only if $q_1 \equiv 0$, $w_1 \equiv 0$ and $\{q_{2|P} = y^2 - w_{2|P} = 0\} \subset \mathbb{P}(P) \cong \mathbb{P}^{10}$ is not a non-singular $8$-dimensional complete intersection for a general $11$-plane $P \subset \mathbb{C}^{M+2}$. To obtain a lower bound for the codimension of $S_{{\rm R2.2}^{\rm 2}}$ in ${\cal F}^o$ we follow the strategy in~\cite[Sec.2.2\&2.3]{EvP2017}; our situation is much simpler but requires some adjustments. \vspace{0.1cm}

{\bf Proposition 2.2.} {\it If $Q = \{q_2=y^2-w_2=0\} \subset \mathbb{P}^{M+1}$ is an irreducible and reduced complete intersction with ${\rm codim}_Q{\rm Sing}(Q) \geqslant 9$ then $Q \cap \mathbb{P}(P)$ is non-singular for a general $11$-dimenional hyperplane $P \subset \mathbb{C}^{M+2}$.} \vspace{0.1cm}

{\bf Proof.} This follows from a version of Bertini's Theorem implying that ${\rm Sing}(Q \cap \mathbb{P}(P)) \subset {\rm Sing}(Q)$ for a general hyperplane $P \subset \mathbb{C}^{M+2}$ (see \cite[II.Thm.8.18]{Hart:AG}), and the fact that a general $10$-dimensional hyperplane $\mathbb{P}(P)$ will not intersect the $\geqslant 11$-codimensional algebraic subset ${\rm Sing}(Q) \subset \mathbb{P}^{M+1}$. Q.E.D.\vspace{0.1cm}

The proposition shows that is enough to find lower bounds for the codimension of the set of pairs $(g,h) \in {\cal F}^o$ such that $q_1 \equiv 0$, $w_1 \equiv 0$ and $Q = \{q_2 = y^2 - w^2 = 0\} \subset \mathbb{P}^{M+1}$ is reducible or non-reduced, and the (Zariski closure of the) set of pairs $(g,h)$ such that $q_1 \equiv 0$, $w_1 \equiv 0$ and ${\rm codim}_Q{\rm Sing}(Q) \leqslant 8$. In both cases we have ${\rm codim}_{\mathbb{P}^{M+1}}Q = 2$ as long as $q_2 \not\equiv 0$ since then $q_2$ cannot have a factor in common with $y^2 - w_2$.

We split up the first set into pairs where the quadric $Q_2 = \{q_2 = 0\} \subset \mathbb{P}^{M+1}$ is reducible or non-reduced, and pairs where $Q_2$ is irreducible and reduced and $Q$ not. $Q_2$ is reducible or non-reduced if and only if $q_2$ is a product of two linear forms. The set of such quadrics has codimension $\binom{M+2}{2} - 2(M+1)$ in ${\cal P}_{2,M+1}$, so the codimension of this component of the first set in ${\cal F}^o$ is
\begin{equation} \label{codimR22-1-1-eq}
2(M+1) + \binom{M+2}{2} - 2(M+1) = \binom{M+2}{2}.
\end{equation}

Next we assume that $Q_2$ is irreducible and reduced. By Grothendiecks Parafactoriality Theorem \cite{CL} and the Lefschetz Theorem for Picard groups \cite[Ex.3.1.35]{Laz1} classes of Weil divisors on $Q_2$ are classes of restrictions of hypersurfaces in $\mathbb{P}^{M+1}$. Furthermore,
$$
H^0(\mathbb{P}^{M+1}, \mathcal{O}_{\mathbb{P}^{M+1}}(a)) \rightarrow H^0(Q_2, \mathcal{O}_{Q_2}(a))
$$
is surjective for all integers $a \geq 0$, bijective for $a \neq 2$ and has kernel $\mathbb{C} \cdot q_2$ for $a=2$. Thus, $Q$ is reducible or non-reduced if and only if $y^2 - w_2 + \lambda q_2$ is a product of linear forms, for some $\lambda \in \mathbb{C}$. But this is only possible if $w_2-\lambda q_2$ is a square of a linear form. For fixed $q_2$ such $w_2$ form a set of codimension $\binom{M+2}{2} - (M+1) - 1$ in ${\cal P}_{2,M+1}$, so the codimension of this component of the first set in ${\cal F}^o$ is
\begin{equation} \label{codimR22-1-2-eq}
2(M+1) + \binom{M+2}{2} - (M+1) - 1 = \frac{(M+4)(M+1)}{2} - 1.
\end{equation}

Now assume that $Q$ is an irreducible and reduced complete intersection of dimension $M-1$ and ${\rm codim}_Q{\rm Sing}(Q) \leqslant 8$. A point $p \in Q$ is a singularity of $Q$ if and only if the tangent space to $Q_2$ in $p$ is contained in the tangent space to $W_2 = \{y^2 - w_2 = 0\}$ in $p$, or vice versa. In both cases there exists a $\lambda = (\lambda_1: \lambda_2) \in \mathbb{P}^1$ such that $W(\lambda) = \lambda_1 q_2 + \lambda_2(y^2-w_2)$ has a singularity in $p$. Thus
$$
{\rm Sing}(Q) \subset \bigcup_{\lambda \in \mathbb{P}^1} {\rm Sing}(W(\lambda)),
$$
and since $\dim {\rm Sing}(Q) \geqslant M-9$ we have $\max_{\lambda \in \mathbb{P}^1} \dim {\rm Sing}(W(\lambda)) \geqslant M-10$. Since $W(\lambda)$ is the vanishing locus of a quadric in $M+2$ variables, $\dim {\rm Sing}(W(\lambda)) = M+1-{\rm rk}(W(\lambda))$, and this implies $\min_{\lambda \in \mathbb{P}^1} {\rm rk}(W(\lambda)) \leqslant 11$.

We distinguish two cases: If ${\rm rk}(q_2) \leq 11$ the inequality above is satisfied for $\lambda = (1:0)$. The codimension of this component of the second set in ${\cal F}^o$ where $q_2$ satisfies this condition is $\geqslant$ to
\begin{equation} \label{codimR22-2-1-eq}
2(M+1) + \frac{(M-9)(M-10)}{2} = \frac{(M-9)(M-6)}{2} + 20.
\end{equation}

If ${\rm rk}(q_2) > 11$ we must find a $\mu \in \mathbb{C}$ such that ${\rm rk}(y^2-w_2 + \mu q_2) \leqslant 11$. This is the case if and only if ${\rm rk}(w_2 - \mu q_2) \leqslant 10$, so for fixed $q_2$ the quadratic polynomial $w_2$ lies in the cone in ${\cal P}_{2,M+1}$ spanned by the vertex $q_2$ and all the quadratic polynomials of rank $\leqslant 10$ in $M+1$ variables. This cone has codimension $\frac{(M-8)(M-9)}{2} - 1$ in ${\cal P}_{2,M+1}$, so the Zariski closure of the set of all pairs $(g,h)$ in ${\cal F}^o$ where $q_2$ and $w_2$ satisfy the above conditions has codimension $\geqslant$ to
\begin{equation} \label{codimR22-2-2-eq}
2(M+1) + \frac{(M-8)(M-9)}{2} - 1 = \frac{(M-8)(M-5)}{2}+17
\end{equation}

\vspace{0.1cm}


{\bf 2.3. Proof of Theorem 2.} Using the estimates (\ref{codimA1-eq}) -- (\ref{codimR22-2-2-eq}) Proposition 2.1 tells us that we will obtain a lower bound for the codimension of the regular locus ${\cal F}_{\rm reg}$ in ${\cal F}$ by subtracting $M$ from the minimum of
$$
\binom{M+1}{2}, \binom{M+2}{2}, \frac{(M-4)(M-3)}{2} + 5,
$$
subtracting $M-1$ from the minimum of
$$
\binom{M+1}{2}, \binom{M+2}{2}, \binom{M}{2}, \frac{(M-4)(M-3)}{2} + 4, \frac{(M-5)(M-4)}{2} + 6, \binom{M+2}{2},
$$
$$
\frac{(M+4)(M+1)}{2} - 1, \frac{(M-9)(M-6)}{2} + 20,  \frac{(M-8)(M-5)}{2} + 17
$$
and taking the smaller of the two numbers. For each $M \geqslant 10$ an elementary calculation yields the lower bound $\xi(M)$ as defined in the Introduction.

\begin{flushleft}
Thomas Eckl and Aleksandr Pukhlikov\\
Department of Mathematical Sciences\\
The University of Liverpool\\
Mathematical Sciences Building\\
Liverpool, L69 7ZL\\
England, U.K
\end{flushleft}

\noindent {\it eckl@liverpool.ac.uk, pukh@liverpool.ac.uk}

\end{document}